\documentclass[twocolumn]{autart}
\usepackage{}

\usepackage{mathrsfs}
\usepackage{amsfonts}
\usepackage{color}                                     

\usepackage[dvips]{graphicx}          

\usepackage{graphicx}

\usepackage{fancyhdr}
\usepackage{amsmath}
\usepackage{amssymb}
\usepackage{mathrsfs}
\usepackage{multicol}
\usepackage{bm}

\usepackage{color}
\usepackage{float}

\renewcommand{\thesection}{\arabic{section}}
\renewcommand{\thesubsection}{\arabic{section}.\arabic{subsection}}

\newtheorem{theorem}{Theorem}[section]
\newtheorem{remark}{Remark}[section]
\newtheorem{lemma}{Lemma}[section]
\newtheorem{example}{Example}[section]
\newtheorem{assumption}{Assumption}[section]
\newtheorem{corollary}{Corollary}[section]
\newtheorem{proposition}{Proposition}[section]
\newtheorem{definition}{Definition}[section]

\emergencystretch=\maxdimen
\usepackage{graphicx}

\begin{document}

\begin{frontmatter}

\title{
Predefined-time Stabilization for Nonlinear  Stochastic Systems\thanksref{footnoteinfo}} 

\thanks[footnoteinfo]{Corresponding author: Shengyuan Xu.\\
}

\author[SCUT]{Tianliang Zhang} \ead{t\_lzhang@163.com},~
\author[SCUT]{Shengyuan  Xu}  \ead{syxu@njust.edu.cn}

\address[SCUT]{School of Automation,  Nanjing University of Science and Technology,
             Nanjing 210094,  Jiangsu Province, China
}



\begin{abstract}                          
\hspace{0.13in}  In this paper, a control scheme for
 stochastic predefined-time stabilization is proposed, which
 improves the control effect compared with  stochastic finite-time or fixed-time stabilization.
The stochastic predefined-time stabilization allows  the upper bound
of the mathematical expectation of  the settling-time function below
an any  given positive value. Some Lyapunov-type results for
predefined-time stabilization  of general stochastic It\^o systems
are presented. Moreover, a state feedback control scheme is designed
for a class of stochastic nonlinear systems in strict-feedback form.
Two simulation examples are supplied to show the usefulness of the
proposed stochastic predefined-time stabilization.

\end{abstract}

\end{frontmatter}

\section{Introduction}

\hspace{0.13in}  The study of convergence time of dynamic  systems is of both practical and theoretical importance. For
some engineering  requirements,  we need to control the convergence
speed faster or slower. Different from  Lyapunov  stability, which is
to study  the system asymptotic  behavior in an infinite time
horizon,
  finite-time or fixed-time stability  investigates the transient response of the system state in a finite time horizon.
  Finite-time control has important applications in robot manipulators \cite{Hongyiguang2002}. Therefore, adaptive finite-time
  tracking control \cite{lihongyi1}, fixed-time stabilization \cite{Polyakov1} and global fast finite-time stabilization of high-order
  nonlinear systems \cite{sunzongyao1}
have been extensively studied in recent years.
 However, in many  practical
applications, it is expected  that the state trajectory  of a
controlled system can converge to the equilibrium point at  any
admissible   time through appropriately  adjusting  control
parameters \cite{Krstic3,Krstic2017,Krstic2019}. A new concept
called predefined-time stability  was proposed in \cite{2015ACC},
which can achieve some consequents that cannot be provided by the
traditional finite-time control schemes, such as the arbitrarily
adjustable upper bound of the convergence time independent of the
initial value. Recently, great progress has been made for
predefined-time control of  deterministic  nonlinear systems
\cite{2020TAC,ND2021,JunkangNi,EJC2020}.

Due to wide applications of stochastic systems,   stochastic
analysis and synthesis  have been  popular research areas over the
past few
 decades
 \cite{Cuixie,hasmin_80,Wuquan,mao,xie2014,xiexuejun1,xiexuejun2013,Yinjuliang1,2015Yin,Yinjuliang2,zhangweihai,ziji1}. Hence,
 it is very valuable to study the finite-time convergence behavior of stochastic systems
 and generalize the above mentioned works to stochastic systems.
Recently, stochastic finite-time stability/stabilization
\cite{Chen2010,wuyuqiang2019,huangshipei,xiexuejun1,Yinjuliang1,2015Yin,Yinjuliang2}
and stochastic fixed-time stability/stabilization
\cite{Hou2021,Songzhibao} have  been studied extensively. Stochastic
finite-time  stability  was first strictly  analyzed in
\cite{Yinjuliang1}, which means that the stochastic settling-time
function $T(x_0)$ is finite almost surely and the system is stable
in probability. Additionally, if we also have $\sup_{x_0\in
{\mathcal R}^n}{\mathcal E}T(x_0)<\infty$, then the stochastic
uncontrolled system is called fixed-time stable. Stochastic
fixed-time stability means that the system  is not only finite-time
stable  but also $\sup_{x_0\in {\mathcal R}^n}{\mathcal
E}T(x_0)<T_{\max}$  for a fixed average time $T_{\max}$ with  any
initial state $x_0\in {\mathcal R}^n$.  In order to obtain a fast
convergence speed, we expect that
 the upper bound of the mathematical expectation of the
settling-time function can  be arbitrarily adjusted by a control
input. Therefore, it is necessary   to  generalize the
predefined-time concept of deterministic systems to stochastic
systems.  However, to the  best of our  knowledge, up to now, there
seems no work  on stochastic predetermined-time stabilization based
on stochastic settling-time function. It should be mentioned that a
recent work \cite{Wuquan} studied  stochastic nonlinear
prescribed-time stabilization in mean square sense.

This paper investigates the problem of stochastic predefined-time stabilization
 of nonlinear stochastic It\^o systems. The main contributions of
this paper are highlighted as follows:
\begin{itemize}
\item [1.] Stochastic predefined-time stability and stabilization are introduced and  a
stochastic predefined-time  stabilization  theorem is obtained for
 general nonlinear stochastic It\^o systems; see
 Theorem~\ref{predefined-time-th1}. Because predefined-time
 stability is stronger than finite-time stability and fixed-time
 stability, Theorem~\ref{predefined-time-th1} is applicable to
 stochastic finite-time/fixed-time stabilization of \cite{Yinjuliang1,Songzhibao}, which can also be viewed as
 extensions of \cite{2018IFAC,2020TAC} to stochastic systems.
Theorem~\ref{predefined-time-th1} is a Lyapunov-type theorem
similarly to the results of \cite{Yinjuliang1,Yinjuliang2}.
Theorem~\ref{predefined-time-th1} yields some useful corollaries
that can be conveniently used in practice.

\item[2.] The high-order nonlinear stochastic system is an important class of stochastic systems that can be used to describe many phenomena
arising from  mechanical systems. Its feedback stabilization has
been investigated in
\cite{Cuixie,Songzhibao,xie2014,xiexuejun1,xiexuejun2013}. Based on
the results of Section 2, we also study the predefined-time
stabilization of high-order nonlinear stochastic systems and give a
practical controller design  algorithm.
  Through  adding the power integrator technique together with a corollary given in Section 2, a
state feedback control scheme is given  for a class of high-order
stochastic nonlinear systems in strict-feedback form, which
guarantees that  the closed-loop system is stochastically
predefined-time stabilizable. In addition, the proposed stochastic predefined-time stabilization control scheme
  improves some existing results \cite{Songzhibao,huangshipei,wuyuqiang2019}.

 \end{itemize}

The rest of this paper is organized as follows:  In Section 2, we make  some preliminaries by introducing some
 new definitions, theorems and corollaries,   Lyapunov-type theorems about the stochastic predetermined-time
 stabilization are obtained.  Section 3 presents the controller  design procedure for high-order nonlinear stochastic
 systems. In Section 4, two simulation examples are given to show the effectiveness of our main results.
 Section 5  concludes this paper with some remarks.

Notation: $\mathcal{R}^n$ denotes the  $n$-dimensional real
Euclidean vector space. $\mathcal{R}_+:=[0,\infty)$. $D'$ is the
transpose of a matrix or vector $D$. $\mathcal{C}^2$ stands for the set of real-valued
twice continuously differentiable functions. $sign(x):=1$ for $x>0$,
$-1$ for $x<0$, $0$ for $x=0$.  For any $b\geq0$, $a\in\mathcal
{R}$, the function $\lceil a\rfloor^b$ is defined as $\lceil
a\rfloor^b=sign(a)|a|^b$. $\mathcal{E}$ means the mathematical
expectation operator. $\mathcal{P}(A)$ represents the probability of
event $A$. $I_{A}(x)$ is the  indicator function, i.e., $I_{A}(x)=1$
for $x\in A$, otherwise, $I_{A}(x)=0$. ${\mathcal K}^1$-functions: A
scalar continuous function $f$ defined from ${\mathcal R}_+$ to
$[0,1)$ is said to be a ${\mathcal K}^1$-function if it is strictly
increasing, $f(0)=0$, and $f(x)\to 1$ as $x\to\infty$.

\section{Stochastic predefined-time stabilization}
In this section, we will consider the following continuous-time stochastic system:
\begin{align}\label{predefined-time-2020-system1}
    &dx(t)=f(x(t),u(t))\,dt+g(x(t),u(t))\,dw(t),\\
    &x(0)=x_0\in \mathcal{R}^n \setminus \{0\},
    \nonumber
\end{align}
where $x(t)\in\mathcal {R}^{n}$ represents the system state.
$u\in\mathcal{R}^{n_u}$ stands for the control input.  $w(t)$ is a
standard one-dimensional Wiener process defined on the filtered
probability space $(\Omega, \mathcal {F}, \{{\mathcal
{F}}_t\}_{t\geq 0}, \mathcal {P})$. $u(t)\in\mathcal {R}^{n_u}$
stands for the control input.  The admissible control
set $\mathcal {U}$  consists of   all ${\mathcal F}_t$-adaptive
control  processes $u(x(t))$, $u(0)=0$,  which makes the
closed-loop system
$$
dx(t)=f(x(t),u(x(t)))\,dt+g(x(t),u(x(t)))\,dw
$$
admits a unique weak solution in forward time for   $x_0\in
\mathcal{R}^n\setminus \{0\}$.  In the considered system, we assume
that
$f:\mathcal{R}^n\times\mathcal{R}^{n_u}\rightarrow\mathcal{R}^n$ and
$g:\mathcal{R}^n\times\mathcal{R}^{n_u}\rightarrow\mathcal{R}^n$ are
continuous in $x$ and $u$ satisfying $f(0,0)=0$ and $g(0,0)=0$. The
purpose of this section is to find an  admissible  control law
$u(x)\in\mathcal {U}$ to stabilize the stochastic system
(\ref{predefined-time-2020-system1}) before a predefined time, i.e.,
achieve the stochastic predefined-time stabilization of system
(\ref{predefined-time-2020-system1}).


\begin{remark}
When studying finite-time stable systems, we are interested in
having a unique solution in forward time \cite{Bhat2000}, which
means that, for any non-zero initial condition $x_0$, $x(t,x_0)$ is
unique before reaching $0$. For stochastic finite-time stable
systems, the concept is extended to the solution in the weak sense
\cite{Yinjuliang2}. Based on the assumption of $f(0,0)=0$ and
$g(0,0)=0$, the origin  is an equilibrium point of
(\ref{predefined-time-2020-system1}). The following lemma gives an
existence result of a weak solution to system
(\ref{predefined-time-2020-system1}). For the trajectory $x(t)$
after reaching the zero equilibrium point, the finite-time
attractiveness (in Definition \ref{predefined-def1qq}) is provided.
\end{remark}

\begin{lemma}\label{predefined-time-2020-newlemma1}\cite{2015Yin,Yinjuliang2}
Suppose that there exits a nonnegative radially unbounded $\mathcal{C}^2$-function $V(x)$ for system
(\ref{predefined-time-2020-system1}), if $\mathcal{L}V(x)\leq0$, then system (\ref{predefined-time-2020-system1}) has a regular continuous solution for any initial value.
\end{lemma}

\begin{remark}
The regular solution means that there is no finite explosion time with probability $1$.
\end{remark}

Before introducing  stochastic predefined-time stabilization,
we will first review some well-known  definitions on stochastic finite-time
(fixed-time, predefined-time, respectively) stabilization.

\begin{definition}\label{predefined-def1qq}\cite{Yinjuliang2,Songzhibao}
System (\ref{predefined-time-2020-system1}) is said to be
stochastically  finite-time stabilizable or finite-time
stabilizable  in probability, if  there exists a state feedback control
$u(t)=u^*(x(t))\in\mathcal {U}$, such that the closed-loop system
\begin{align}\label{predefined-time-2020-closed system1}
    &dx(t)=f(x(t),u^*(x(t)))\,dt+g(x(t),u^*(x(t)))\,dw(t),\\
    &x(0)=x_0\in \mathcal{R}^n\setminus \{0\}, \nonumber
\end{align}
is stochastically finite-time stable, i.e.,
\begin{itemize}
  \item Finite-time attractiveness: For any non-zero initial value $x_0$, there exists the first settling-time function
    $T(x_0):=\inf\{t: x(t,u^*,x_0)=0\}$, such that
    $$
    \mathcal{P}(T(x_0)<\infty)=1,
    $$
where  $x(t,u^*,x_0)$ is the solution of
(\ref{predefined-time-2020-closed system1}). Moreover, $x(t+T(x_0),u^*,x_0)\equiv0$, a.s., for any $t\geq0$.

  \item Stability  in probability: For any pair $(\alpha, \beta)$, $\alpha\in(0,1)$, $\beta>0$,
  there exists a $\delta(\alpha, \beta)>0$, such that  for  all $\|x_0\|<\delta(\alpha, \beta)$, we have
  \begin{align*}
      \mathcal{P}(\|x(t, u^*, x_0)\|<\beta,  \forall t>0)\geq1-\alpha.
  \end{align*}
  Moreover,  if there exists a positive constant $T_{max}$ such  that
\begin{align*}
    \sup_{x_0\in\mathcal{R}^n\setminus \{0\}}\mathcal{E}T(x_0)\leq T_{max},
\end{align*}
then  system (\ref{predefined-time-2020-closed system1}) is said to be
stochastically fixed-time stable, and  system
(\ref{predefined-time-2020-system1}) is said to be stochastically fixed-time
stabilizable.
\end{itemize}
\end{definition}

\begin{remark}
In recent literature, new and improved definitions called
predefined-time stable and stabilizable were  proposed for
deterministic systems based on  fixed-time stability
\cite{2020TAC,2015ACC}. This paper aims  to extend these  concepts
and related results  to stochastic system
(\ref{predefined-time-2020-system1}).
\end{remark}

\begin{definition}\label{predefined-def2}
System (\ref{predefined-time-2020-system1}) is said to be stochastically predefined-time stabilizable
if it is stochastically fixed-time stabilizable and  for any $\varsigma>0$, there exists a control
$u(t)=u^*(x(t))\in\mathcal {U}$, such that
\begin{align*}
    \sup_{x_0\in\mathcal{R}^n\setminus \{0\}}\mathcal{E}T(x_0)\leq \varsigma.
\end{align*}
\end{definition}

\begin{remark}
Recently,   another newly proposed definition  called
prescribed-time mean-square stability  was introduced \cite{Wuquan},
which  is different from Definition~\ref{predefined-def2}.
Definition~\ref{predefined-def2} is based on stochastic settling
time function, when the system degenerates into deterministic
systems,  Definition~\ref{predefined-def2}  is   consistent  with
the corresponding definition of deterministic systems.
\end{remark}

\begin{remark}
From Definitions \ref{predefined-def1qq} and \ref{predefined-def2},
for the  system (\ref{predefined-time-2020-system1}) with non-zero
initial value, there must be $x(t+T(x_0), u^*, x_0)\equiv0$ with
$t\geq0$ and $T(x_0)$ being the stochastic settling-time function.
So, $0$ is an absorbing state. Moreover, when discussing
predefined-time stabilization or finite-time stabilization control
problems, a non-zero initial value is usually assumed
\cite{2020TAC,2015Yin}.
\end{remark}

The following theorem is a Lyapunov-type theorem about
stochastic predefined-time stabilization.

\begin{theorem}\label{predefined-time-th1}
For any $\alpha>0$, if there exists a control input $u^*(x(t))\in\mathcal {U}$,
driving the state of system (\ref{predefined-time-2020-system1}) to
satisfy
\begin{equation}\label{predefined-time-con1}
\mathcal{L}_{u=u^*}V(x)\leq-\frac{1}{\alpha}\beta(V(x)), \
 x\in\mathcal{R}^n \setminus \{0\},
\end{equation}
where $\beta:\mathcal{R}_{+}\mapsto\mathcal{R}_{+}$ with
$\dot{\beta}(\cdot)\geq0$, $\beta(s)>0$ for any $s>0$, and
$\int^{\infty}_0\frac{1}{\beta(s)}\,ds\leq1$,
$V:\mathcal{R}^n\mapsto\mathcal{R}_{+}$ is a
$\mathcal{C}^2$-positive definite and radially unbounded function, and $\mathcal{L}_{u=u^*}$ represents the infinitesimal generator of system (\ref{predefined-time-2020-closed system1}).
Then system (\ref{predefined-time-2020-system1}) is stochastically
predefined-time stabilizable, and $\sup_{x_0\in\mathcal{R}^n\setminus \{0\}}\mathcal{E}T(x_0)\leq \alpha.
$.
\end{theorem}
{\textbf {Proof.}}  By Theorem 1 of \cite{Yinjuliang2},  system
(\ref{predefined-time-2020-system1}) is stochastically finite-time
stabilizable under the conditions of this theorem. Hence, in order
to prove stochastic predefined-time stabilization of system
(\ref{predefined-time-2020-system1}), we only need  to  show that
for any $\alpha>0$, the following holds:
$$
\sup_{x_0\in\mathcal{R}^n\setminus \{0\}}\mathcal{E}T(x_0)\leq \alpha.
$$

Obviously,  if $x_0=0$, it directly leads to  $x(t,u^*,x_0)=0$ a.s.
for any $t\geq0$. So we only need to consider the case of non-zero
initial state. From condition
(\ref{predefined-time-con1}) and Lemma
\ref{predefined-time-2020-newlemma1}, it leads to that for each
$x_0\neq0$, there exists a regular continuous solution $x(t,u^*,
x_0)$ to (\ref{predefined-time-2020-system1}). Therefore, there
exists a minimal positive integer $k_0\in \{1,2,\cdots,\}$  such
that $\frac{1}{k_0}<\|x_0\|<k_0$. Define some sequences as follows:
\begin{align*}
    &\tau_{k}=\inf\left\{t\geq0:\|x(t,u^*, x_0)\|\notin \left(\frac{1}{k},k\right)\right\},\\
    &\tau_{1k}=\inf\left\{t\geq0:\|x(t,u^*,x_0)\| \in\left[0,\frac{1}{k}\right]\right\},\\
    &\tau_{2k}=\inf\left\{t\geq0:\|x(t, u^*, x_0)\| \in\left[k,\infty\right)\right\}.
\end{align*}
 Considering Lemma \ref{predefined-time-2020-newlemma1}
and the definition of the  admissible control set $\mathcal {U}$,
according to the definition of stopping time, sequences $\tau_{k}$,
$\tau_{1k}$ and  $\tau_{2k}$ are  ${\mathcal F}_{\tau_{k}}$-,
${\mathcal F}_{\tau_{1k}}$-  and  ${\mathcal
F}_{\tau_{2k}}$-measurable, respectively. So they are all  stopping
times. For convenience, in the sequel, we denote the solution
$x(t,u^*,x_0)$ by $x(t)$ for short. We introduce a new Lyapunov
function $P(x)=\int^{V(x)}_0\frac{\alpha}{\beta(s)}ds$.  By It\^o
formula, we get
\begin{align*}
dP(x(t))&=\mathcal{L}_{u=u^*}P(x(t))\,
dt\\
 &\ \ +\frac{\alpha}{\beta(V(x(t)))}\frac{\partial V'}{\partial
x}g(x(t),u^*(x(t)))\,dw(t),
\end{align*}
where $\mathcal{L}_{u=u^*}P(x)$ can be computed as
\begin{align}\label{predefined-time-LP}
  &\mathcal{L}_{u=u^*} P(x)=\frac{\alpha\mathcal{L}_{u=u^*} V(x)}{\beta(V(x))}\nonumber\\
  &-\frac 1 2
  \frac{\alpha\dot{\beta}(V)}{\beta^2(V)}\left(g'(x,u^*)\frac{\partial V}{\partial x}\right)\left(\frac{\partial V'}{\partial x}g(x,u^*)\right).
\end{align}
Then, similar to Theorem 3.1 in \cite{Yinjuliang1}, it follows  that
\begin{align}
    \mathcal{E}P(x(t\wedge\tau_k))-P(x_0)=\mathcal{E}\int^{t\wedge\tau_k}_0\mathcal{L}_{u=u^*}P(x(v))dv.
\end{align}
From (\ref{predefined-time-con1}) and (\ref{predefined-time-LP}), we
have $\mathcal{L}_{u=u^*} P(x)\leq-1$ for any $x\neq0$. So
$\mathcal{E}P(x(t\wedge\tau_k))-P(x_0)\leq-\mathcal{E}(t\wedge\tau_k)$ and
$\mathcal{E}(t\wedge \tau_k)\leq P(x_0)$ hold.  Because $\{t\wedge
\tau_k\}_{k\in \{k_0,k_0+1,\cdots\}}$ is an increasing sequence of
nonnegative random variables,
 by monotonic convergence theorem,
\begin{align*}
 \lim_{k\rightarrow\infty}\mathcal{E}(t\wedge \tau_k)=\mathcal{E}(\lim_{k\rightarrow\infty} t\wedge \tau_k)=\mathcal{E}(t\wedge
\tau_{\infty}).
\end{align*}
Note that
$$ \tau_{\infty}=\tau_{1,\infty}\wedge\tau_{2,\infty}.
$$
Since $x(t,x_0)$ is a  regular solution, $\tau_{2,\infty}=\infty$
a.s..
So $\mathcal{E}(t\wedge
\tau_{\infty})=\mathcal{E}(t\wedge\tau_{1,\infty})\leq
P(x_0)\le \alpha$. Due to the arbitrariness of $x_0$ and $\tau_{1,\infty}=T(x_0)$, it yields
that
$$
\sup_{x_0\in\mathcal{R}^n\setminus \{0\}}\mathcal{E}T(x_0)\leq \alpha.
$$
From (\ref{predefined-time-con1}) and $V(x)\geq0$, $V_t:=V(x(t,u^*,x_0))$ is a nonnegative continuous supermartingale with augmented filtration $\{{\mathcal {F}}_t\}_{t\geq 0}$. Through Doob's optional-sampling theorem for continuous nonnegative supermartingales in \cite{Rogers2000},
\begin{align*}
\mathcal {E}(V_{T(x_0)+t}|\mathcal {F}_{T(x_0)})\leq V_{T(x_0)}=0, \forall t\geq0.
\end{align*}
Taking expectation on both sides of the above inequality, we have $\mathcal {E}V_{T(x_0)+t}\equiv0$. Since $V$ is positive definite, it follows that $x(T(x_0)+t,u^*,x_0)\equiv0$, a.s..
The
proof is completed.

Generally speaking, stochastic predefined-time stabilization is a
special case of stochastic fixed-time stabilization, which is
stronger than stochastic fixed-time stabilization. In order to
illustrate the difference between these two concepts, a numerical
example is presented below.

\begin{example}
Consider a scalar system
\begin{equation}\label{predefined-ee1}
    dx=udt+xdw, \  x(0)=x_0.
\end{equation}
 By selecting the control input
\begin{equation}\label{predefined-uu1}
 u=u^*=-\frac{1}{2}x-x^{a}-x^b
\end{equation}
with $a\in(0,1)$ and $b>1$, then  the  feedback system of
(\ref{predefined-ee1}) becomes
\begin{equation}\label{predefined-ee1feed}
    dx=(-\frac{1}{2}x-x^{a}-x^b)\,dt+x\,dw, \  x(0)=x_0.
\end{equation}
For system  (\ref{predefined-ee1feed}) and the Lyapunov function
$V(x)=x^2$,   it is easy to compute
$$
{\mathcal L}_{u=u^*}V(x)=-2x^{a+1}-2x^{b+1},
$$
where ${\mathcal L}_{u=u^*}$ is the infinitesimal generator of
system  (\ref{predefined-ee1feed}).  According to Lemma 7 of
\cite{Hou2021}, we have
\begin{align*}
    \sup_{x_0\in\mathcal{R}}\mathcal{E}T(x_0)\leq T_{max}:=\frac{1}{1-a}+\frac{1}{b-1}.
\end{align*}
Hence,  system (\ref{predefined-ee1}) is  stochastically fixed-time
stabilizable.

In addition, for any $\alpha>0$,  if we choose the control input as
\begin{equation}\label{predefined-uu2}
u=u^*(x)=-\frac {\sqrt\pi} {2\alpha}sign(x)e^{x^2}-\frac 1 2x,
\end{equation}
and $\mathcal{C}^2$-positive definite and radially unbounded
function $V(x)=x^2$, then by  It\^o formula, we have
\begin{align*}
  \mathcal{L}_{u=u^*}V(x)&=2xu^*+x^2=x(-\frac {\sqrt\pi}
  {\alpha}sign(x)e^{x^2})\\
  &=-\frac {\sqrt\pi}
  {\alpha}|x|e^{x^2}=-\frac {\sqrt\pi}{\alpha} V(x)^{1/2}e^{V(x)}.
\end{align*}
Let $\beta(s)=\sqrt \pi s^{1/2}e^s$,  then $\beta(s)$ satisfies the
conditions of Theorem~\ref{predefined-time-th1}. By
Theorem~\ref{predefined-time-th1}, system (\ref{predefined-ee1}) is
also stochastically predefined-time stabilizable, and
(\ref{predefined-uu2}) is a stochastic predefined-time stabilizing
controllor.
\end{example}

\begin{remark}
Fixed-time stabilization is often difficult and sometimes impossible to adjust the controller gain such that achieving stabilization within a required predefined time. There is no guarantee that this upper bound can be adjusted arbitrarily. On this point, we can refer to the discussion in \cite{2020TAC,EJC2020}. From the above example, it is easier to tune the convergence time bound for stochastic predefined-time stabilization than stochastic fixed-time stabilization. Therefore, the stochastic fixed-time control discussed in \cite{Songzhibao} cannot achieve the effect of the stochastic predefined-time control discussed in this paper. Furthermore, this paper proves various Lyapunov-type theorems for stochastic predefined-time stabilization, and presents a design scheme of stochastic predefined-time controllers for higher-order nonlinear systems.
\end{remark}

\begin{remark}\label{predefined-time-2020-rem1}
Theorem \ref{predefined-time-th1} presents a Lyapunov-type result
about stochastic predefined-time stabilization. When
$\int^{t}_0\frac{1}{\beta(s)}ds<\infty$ with $t\geq0$ replaces
$\int^{\infty}_0\frac{1}{\beta(s)}ds\leq1$ in Theorem
\ref{predefined-time-th1}, then system
(\ref{predefined-time-2020-system1}) is stochastically finite-time
stabilizable \cite{Yinjuliang1}. In addition, if the condition
$\int^{\infty}_0\frac{1}{\beta(s)}ds\leq 1$ is changed as
$\int^{\infty}_0\frac{1}{\beta(s)}\,ds<\infty$, then system
(\ref{predefined-time-2020-system1}) is stochastically fixed-time
stabilizable.
\end{remark}

\begin{corollary}\label{predefined-time-cor1}
For any $\alpha>0$, if there exists a control input $u^*(x(t))\in \mathcal {U}$, driving the state of system (\ref{predefined-time-2020-system1}) to
satisfy
\begin{equation}\label{predefined-time-con2}
  \mathcal{L}_{u=u^*}V(x)\leq-\frac{\tilde\beta(V(x))^p}{(1-p)\alpha \dot{\tilde \beta}(V(x))},  x\in\mathcal{R}^n \setminus \{0\},
\end{equation}
where $\tilde \beta(\cdot)$ is a function belonging  to
$\mathcal{K}^1$ with $\ddot{\tilde \beta}(\cdot)\leq0$, $0\leq p<1$.
$V:\mathcal{R}^n\mapsto\mathcal{R}_{+}$ is  a \
$\mathcal{C}^2$-positive definite and radially unbounded function.
Then system (\ref{predefined-time-2020-system1}) is stochastically
predefined-time stabilizable,  and for any $x_0\in\mathcal{R}^n\setminus\{0\}$,
the first settling-time function satisfies $\mathcal{E}T(x_0)\leq
\alpha$.
\end{corollary}

{\textbf {Proof.}} Corollary \ref{predefined-time-cor1} can be
directly obtained by Theorem \ref{predefined-time-th1} via setting
$$
\beta(s)=\frac{\tilde \beta(s)^p}{(1-p) \dot{\tilde \beta}(s)}.
$$

We can also  give another Lyapunov-type theorem  inspired by
\cite{2018IFAC}.

\begin{corollary}
\label{predefined-time-cor2}
For any $\alpha>0$, if there exists a control input $u^*(x(t))\in \mathcal {U}$, driving the state of system (\ref{predefined-time-2020-system1}) to
satisfy
\begin{equation}\label{predefined-time-con3}
 \mathcal{L}_{u=u^*} W(x)\leq-\frac{1}{\alpha},\ x\in\mathcal{R}^n \setminus \{0\},
\end{equation}
where $W:\mathcal{R}^n\mapsto \mathcal{R}_{+}$ is a
$\mathcal{C}^2$-positive definite function with $0\leq W(s)<1$, $\lim_{s\rightarrow
\infty}W(s)=1$. Then system (\ref{predefined-time-2020-system1}) is stochastically
predefined-time stabilizable,  and for any $x_0\in\mathcal{R}^n\setminus\{0\}$,
the first settling-time function satisfies $\mathcal{E}T(x_0)\leq
\alpha$.
\end{corollary}

{\textbf {Proof.}} Set  $V(x)=-\ln(1-W(x))$. By
(\ref{predefined-time-con3}), it is easy to obtain that
\begin{equation}\label{eqnjnrfkj}
\mathcal{L}_{u=u^*}W(x)=e^{-V(x)}\mathcal{L}_{u=u^*}V(x)\leq-\frac{1}{\alpha}.
\end{equation}
This gives
$$
\mathcal{L}_{u=u^*}V(x)\le -\frac{1}{\alpha} e^{V(x)}.
$$
By setting $\beta(s)=e^s$, and using Theorem \ref{predefined-time-th1}, we have the desired results
immediately.

\begin{corollary}
\label{predefined-time-cor3} Assume there exist a
$\mathcal{C}^2$-positive definite and radially unbounded function
$V$, positive constants $a>0$, $0<b_1<1$, and $b_2>1$ satisfying
$a-ab_1-1>0$. If for any $\alpha>0$, there exists a control
input $u^*(x(t))\in\mathcal {U}$, such that
\begin{align}\label{predefined-time-con4}
\mathcal{L}_{u=u^*}V(x)
\leq&-\frac{a-ab_1}{\alpha
 (a-ab_1-1)(b_2-1)}V(x)^{b_2}\nonumber\\
 &-\frac{a}{\alpha }V(x)^{b_1}
\end{align}
for $x\in\mathcal{R}^n \setminus \{0\}$, then   system (\ref{predefined-time-2020-system1}) is stochastically
predefined-time stabilizable, i.e.,  $\sup_{\forall
x_0}\mathcal{E}T(x_0)\leq \alpha$.
\end{corollary}
{\textbf {Proof.}} In Theorem \ref{predefined-time-th1}, set
 $$
\beta(V(x))=aV(x)^{b_1}+\frac{a-ab_1}{(a-ab_1-1)(b_2-1)}V(x)^{b_2}.
$$
Obviously, $\beta(s)\ge0$ for $s\in {\mathcal R}_+$, and
$\dot{\beta}(\cdot)\ge 0$.  Moreover,
\begin{align*}
&\int^{\infty}_{0}\frac{1}{\beta(s)}ds=\int^{1}_{0}\frac{1}{\beta(s)}ds+\int^{\infty}_{1}\frac{1}{\beta(s)}\,ds\\
\leq& \int^{1}_{0} \frac {1} {a s^{b_1}}\,ds+\int^{\infty}_{1} \frac{(a-ab_1-1)(b_2-1)}{(a-ab_1)s^{b_2}}\,ds\\
=&\frac{1}{(1-b_1)a}+\frac{a-ab_1-1}{a-ab_1}=1.
\end{align*}
Therefore, (\ref{predefined-time-con4}) implies that conditions of
Theorem \ref{predefined-time-th1} are satisfied. This corollary is
proved.

\section{Controller design in strict-feedback form}
In this section, we address the stochastic predefined-time stabilization
problem for nonlinear stochastic systems based on backstepping method. Consider the following high-order stochastic nonlinear system described  by
\begin{align}\label{predefined-time-2020-system2}
    \begin{cases}
    dx_i=(h_i\lceil x_{i+1}\rfloor^{q_i}+f_i(\bar{x}_i))dt+g_i(\bar{x}_i)dw, \\
    dx_n=(h_n\lceil u\rfloor^{q_n}+f_n(x))dt+g_n(x)dw,
    \end{cases}
\end{align}
where $i=1,2,\cdots,n-1$, $x=[x_1,\cdots,x_n]'\in \mathcal{R}^n$ stands for the system state, $u\in\mathcal{R}$ is the  control input. $w$ is a one-dimensional  standard Wiener process   defined on the filtered probability space
 $(\Omega, {\mathcal F}$, $\{{\mathcal F}_t\}_{t\ge 0}$, ${\mathcal P})$. $\bar{x}_i$ is defined as $\bar{x}_i=[x_1,\cdots,x_i]'$. $h_1,\cdots,h_n$ are unknown virtual control coefficients. The drift terms $f_i(x)$ and diffusion terms $g_i(x)$, $i=1,2,\cdots,n$, are Borel measurable continuous functions with $f_i(0)=0$ and $g_i(0)=0$. For any $i$,
 $1<q_i\in\mathcal{R}_+$ is called the high-order of system (\ref{predefined-time-2020-system2}).

The following  assumption allows  the functions $f_i(\bar{x}_i)$ and $g_i(\bar{x}_i)$ to have high-order and low-order nonlinear growth rates.

\begin{assumption}\label{predefined-time-2020-ass1}
For any $i=1,2,\cdots,n$, the drift term   $f_i(\bar{x}_i)$ and the diffusion term $g_i(\bar{x}_i)$ satisfy the following conditions:
\begin{align*}
   |f_i(\bar{x}_i)|\leq \phi_i(\bar{x}_i)\sum^i_{j=1}|x_j|^{\varpi_{ij}+\frac{r_i+\kappa}{r_j}},\\
   |g_i(\bar{x}_i)|\leq \varphi_i(\bar{x}_i)\sum^i_{j=1}|x_j|^{\rho_{ij}+\frac{2r_i+\kappa}{2r_j}},
\end{align*}
where $\varpi_{ij}\geq 0$, $\rho_{ij}\geq 0$, $\phi_i(\bar{x}_i)$ and $\varphi _i(\bar{x}_i)$ are non-negative smooth functions with $\phi_i(0)=0$ and $\varphi _i(0)=0$,  $\kappa\in(-\frac{1}{1+\sum^{n-1}_{s=1}q_1\cdots q_s},0)$.
$r_{i+1}$ is recursively defined as $r_{i+1}=\frac{r_i+\kappa}{q_{i}}$ with $r_1=1$.
\end{assumption}

\begin{remark}
Assumption \ref{predefined-time-2020-ass1} is a common nonlinear growth rate \cite{wuyuqiang2019,sunzongyao1}. the powers in growth condition of $f_i(\bar{x}_i)$ and $g_i(\bar{x}_i)$ are defined as $\varpi_{ij}+\frac{r_i+\kappa}{r_j}\in\mathcal{R}_+$ and $\rho_{ij}+\frac{2r_i+\kappa}{2r_j}\in\mathcal{R}_+$, respectively. $\varpi_{ij}+\frac{r_i+\kappa}{r_j}$ and $\rho_{ij}+\frac{2r_i+\kappa}{2r_j}$  belonging to an interval $(0,\infty)$ allow $f_i(\bar{x}_i)$ and $g_i(\bar{x}_i)$ to have both both high-order and low-order nonlinear growth rates. Assumption \ref{predefined-time-2020-ass1} includes the counterparts in the closely
related works as special cases.
For example, when $\varpi_{ij}=\rho_{ij}=0$, Assumption \ref{predefined-time-2020-ass1} degenerates to Assumption 1 in \cite{Fang2019}.
If  $\phi_i(\bar{x}_i)$ and $\varphi_i(\bar{x}_i)$ are further specialised as constants, Assumption \ref{predefined-time-2020-ass1} becomes the low-order growth rate used in \cite{xiexuejun1}. When $\varpi_{ij}=1-\frac{r_i+\kappa}{r_j}$ and $\rho_{ij}=1-\frac{2r_i+\kappa}{2r_j}$, Assumption \ref{predefined-time-2020-ass1}  reduces the linear-like growth rate used in \cite{Khoo2013}.
\end{remark}

\begin{assumption}\label{predefined-time-2020-ass2}
For any $i=1,2,\cdots,n$, there exist constants $\overline{h}_i>0$ and $\underline{h}_i>0$ such that $\underline{h}_i\leq h_i\leq \overline{h}_i$.
\end{assumption}

\begin{lemma}\label{predefined-time-2020-lem1}\cite{xie2014}
For any $a\in\mathcal{R}$ and $b\geq2$,  we have that the function $f(a):=\lceil a\rfloor^b\in\mathcal{C}^2$, $\frac{\partial f(a)}{\partial a}=b|a|^{b-1}$ and $\frac{\partial^2 f(a)}{\partial a^2}=b(b-1)\lceil a\rfloor^{b-2}$.
\end{lemma}

\begin{lemma}\label{predefined-time-2020-lem2}\cite{Ding2012}
Suppose  $0<p<1$ and $q>1$, then  for any $x,y\in\mathcal{R}$, we have
\begin{align*}
    |\lceil x\rfloor^{pq}-\lceil y\rfloor^{pq}|\leq2^{1-p}|\lceil x\rfloor^{q}-\lceil y\rfloor^{q}|^p.
\end{align*}
\end{lemma}

\begin{lemma}\label{predefined-time-2020-lem3}\cite{Qian2001}
For any positive real numbers $p$, $q$ and any real-valued function $f(x,y)>0$,  the  following  relationship  holds:
\begin{align*}
    |x|^p|y|^q\leq\frac{p}{p+q}f(x,y)|x|^{p+q}+\frac{q}{p+q}f^{-\frac{p}{q}}(x,y)|y|^{p+q}.
\end{align*}
\end{lemma}

\begin{lemma}\label{predefined-time-2020-lem4}\cite{Qian2001}
For any $x,y\in\mathcal{R}$, $a\geq1$, we have
$(|x|+|y|)^{\frac{1}{a}}\leq|x|^{\frac{1}{a}}+|y|^{\frac{1}{a}}\leq2^{\frac{a-1}{a}}(|x|+|y|)^{\frac{1}{a}}$.
\end{lemma}

\begin{lemma}\label{predefined-time-2020-lem5}\cite{Qian2001}
$(\sum^j_{i=1}a_i)^b\leq\max\{j^{b-1},1\}(\sum^j_{i=1}a_i^b)$ always holds for any positive numbers $a_1, a_2, \cdots, a_j$ and $b$.
\end{lemma}

Now, we are ready to present our  main theorem of this section.
\begin{theorem}\label{predefined-time-th2}
If the high-order stochastic nonlinear system (\ref{predefined-time-2020-system2}) satisfies Assumptions \ref{predefined-time-2020-ass1} and \ref{predefined-time-2020-ass2}, then there exists a state-feedback controller such that the closed-loop system is stochastically predefined-time stabilizable.
\end{theorem}
{\textbf {Proof.}}
The proof is based on  inductive arguments.
{\bf Step 1:}
Firstly, by the power integrator technique, we choose the Lyapunov function for the
 first subsystem as $V_1(x_1)=\int^{x_1}_{x_1^*}\lceil\lceil s \rfloor^{\frac{r}{r_1}}-\lceil x_1^*\rfloor^{\frac{r}{r_1}}
 \rfloor^{\frac{4 r-\kappa-r_1}{r}}\,ds$, where $r\geq \max_{1\leq i\leq n}\{2r_i\}$ and $x^*_1=0$.
 By It\^o's  formula, Assumption \ref{predefined-time-2020-ass1}, and  Lemma \ref{predefined-time-2020-lem1},
 we can get that
 {\small
\begin{align*}
    &\mathcal{L}_1V_1(x_1)\\
    =&\frac{\partial V_1}{\partial x_1}(h_1\lceil x_2\rfloor^{q_1}+f_1(\bar{x}_1))+\frac{1}{2}\frac{\partial^2 V_1}{\partial x_1^2}|g_1(\bar{x}_1)|^2\\
    =&\lceil\xi_1\rfloor^{\frac{4 r-\kappa-r_1}{r}}(h_1\lceil x_2\rfloor^{q_1}+f_1(\bar{x}_1))+\frac{4 r-\kappa-r_1}{2r}|\xi_1|^{\frac{3r-\kappa-r_1}{r}}\\
    &\frac{r}{r_1}|x_1|^{\frac{r}{r_1}-1}|g_1(\bar{x}_1)|^2\\
    \leq&\lceil\xi_1\rfloor^{\frac{4 r-\kappa-r_1}{r}}h_1(\lceil x_2\rfloor^{q_1}-\lceil x_2^*\rfloor^{q_1})+\lceil\xi_1\rfloor^{\frac{4 r-\kappa-r_1}{r}}h_1\lceil x_2^*\rfloor^{q_1}\\
    &+\frac{4 r-\kappa-r_1}{2r}|\xi_1|^{\frac{3r-\kappa-r_1}{r}}\frac{r}{r_1}|x_1|^{\frac{r}{r_1}-1}(\bar{\varphi} _1(\bar{x}_1)|x_1|^{\frac{2r_1+\kappa}{2r_1}})^2\\
    &+\lceil\xi_1\rfloor^{\frac{4 r-\kappa-r_1}{r}}\bar{\phi}_1(x_1)|x_1|^{\frac{r_1+\kappa}{r_1}}\\
    \leq&\lceil\xi_1\rfloor^{\frac{4 r-\kappa-r_1}{r}}h_1(\lceil x_2\rfloor^{q_1}-\lceil x_2^*\rfloor^{q_1})+\lceil\xi_1\rfloor^{\frac{4 r-\kappa-r_1}{r}}h_1\lceil x_2^*\rfloor^{q_1}\\
    &+|\xi_1|^{4}\bar{\phi}_1(x_1)+\frac{4 r-\kappa-r_1}{2r_1}|\xi_1|^{4 }\bar{\varphi}_1^2(\bar{x}_1),
\end{align*}
} where $\xi_1=\lceil x_1\rfloor^{\frac{r}{r_1}}$,
$\bar{\varphi}_1(x_1)\geq\varphi_1(x_1)|x_1|^{\rho_{11}}$ and
$\bar{\phi}_1(x_1)\geq\phi_1(x_1)|x_1|^{\varpi_{11}}$ are smooth
functions. $\mathcal{L}_1$ is the infinitesimal generator of the
first subsystem. The definitions of $\bar{\varphi}_i(\bar{x}_i)$ and
$\bar{\phi}_i(\bar{x}_i)$ used in the later proof are similar and
will not be repeated. Designing
$x_{2}^*=-\beta_1\lceil\xi_1\rfloor^{\frac{r_2}{r}}$ and
$\beta_1>0$, it leads to
\begin{align*}
\lceil x_2^*\rfloor^{q_1}=&sign(-\beta_1\lceil\xi_1\rfloor^{\frac{r_2}{r}})|\beta_1\lceil\xi_1\rfloor^{\frac{r_2}{r}}|^{q_1}\\
=&-sign(\lceil\xi_1\rfloor^{\frac{r_2}{r}})|\beta_1\lceil\xi_1\rfloor^{\frac{r_2}{r}}|^{q_1}\\
=&-\beta_1^{q_1}\lceil\xi_1\rfloor^{\frac{r_2q_1}{r}}.
\end{align*}
Note that $\lceil\xi_1\rfloor^{\frac{4 r-\kappa-r_1}{r}}\lceil\xi_1\rfloor^{\frac{r_2q_1}{r}}=|\xi_1|^4$.
Set
$
\beta_1(x_1)=\left(\frac{4r-\kappa-r_1}{2r_1\overline{h}_1}\bar{\varphi}_1^2(x_1)
+\frac{nk_1}{\overline{h}_1k_4}+\frac{k_2}{\overline{h}_1k_4}|\xi_1|^{k_3}
+\frac{\bar{\phi}_1(x_1)}{\overline{h}_1}\right)^{\frac{1}{q_1}}
$
with $k_1$, $k_2$, $k_3$, $k_4$  to be designed in the future. Then
\begin{align*}
\mathcal{L}_1V_1(x_1)\leq&\lceil\xi_1\rfloor^{\frac{4 r-\kappa-r_1}{r}}h_1(\lceil x_2\rfloor^{q_1}-\lceil x_2^*\rfloor^{q_1})\\
&-\frac{nk_1}{k_4}|\xi_1|^4-\frac{k_2}{k_4}|\xi_1|^{4+k_3}.
\end{align*}
{\bf Step 2-Inductive assumption:}  Suppose for $j\in\{3,4,\cdots,n\}$,  there is a $\mathcal{C}^2$- positive definite radially unbounded function $V_{j-1}$, and a set of virtual controllers $x_1^*, x_2^*,\cdots, x_j^*$ defined by
\begin{align*}
    \begin{array}{cccc}
        x_1^*=0,& \xi_1=\lceil x_1\rfloor^{\frac{r}{r_1}}-\lceil x_1^*\rfloor^{\frac{r}{r_1}};\\
        x_{2}^*=-\beta_1(\bar{x}_1)\lceil\xi_1\rfloor^{\frac{r_2}{r}},& \xi_2=\lceil x_2\rfloor^{\frac{r}{r_2}}-\lceil x_2^{*}\rfloor^{\frac{r}{r_2}};\\
        \vdots &\vdots \\
        x_{j}^*=-\beta_{j-1}(\bar{x}_{j-1})\lceil\xi_{j-1}\rfloor^{\frac{r_j}{r}},& \xi_j=\lceil x_j\rfloor^{\frac{r}{r_j}}-\lceil x_j^{*}\rfloor^{\frac{r}{r_j}}
    \end{array}
\end{align*}
with functions $\beta_i(\bar{x}_i)>0$, $i\in\{1,2,\cdots,j-1\}$, such that
{\small
\begin{align}\label{eqbjbj}
\mathcal{L}_{j-1}V_{j-1}(\bar{x}_{j-1})\leq&\lceil\xi_{j-1}\rfloor^{\frac{4 r-\kappa-r_{j-1}}{r}}h_{j-1}(\lceil x_j\rfloor^{q_{j-1}}-\lceil x_j^*\rfloor^{q_{j-1}})\nonumber\\
&-(n-j+2)\sum^{j-1}_{i=1}\frac{k_1}{k_4}|\xi_i|^4-\sum^{j-1}_{i=1}\frac{k_2}{k_4}|\xi_i|^{4+k_3},
\end{align}}
where $\mathcal{L}_{j-1}$ is the infinitesimal generator of the first $j-1$ subsystems.

{\bf Step 3:} In the following, we shall show that (\ref{eqbjbj}) still holds when $j-1$ is replaced by $j$ with
\begin{align}\label{predefined-time-2020-vv1}
  V_j(\bar{x}_j)=V_{j-1}(\bar{x}_{j-1})+W_j(\bar{x}_j)
\end{align}
and
\begin{align}\label{predefined-time-2020-vv2}
W_j(\bar{x}_j)=\int^{x_j}_{x_j^*}\lceil\lceil s \rfloor^{\frac{r}{r_j}}-\lceil x_j^{*}\rfloor^{\frac{r}{r_j}}\rfloor^{\frac{4 r-\kappa-r_j}{r}}\,ds.
\end{align}
With the help of Proposition \ref{predefined-time-2020-pro1} in Appendix, it can be deduced from $V_j(\bar{x}_j)$ that
\begin{align}\label{predefined-time-2020-2}
&\mathcal{L}_jV_{j}(\bar{x}_{j})\nonumber\\
\leq&\lceil\xi_{j-1}\rfloor^{\frac{4 r-\kappa-r_{j-1}}{r}}h_{j-1}(\lceil x_j\rfloor^{q_{j-1}}-\lceil x_j^*\rfloor^{q_{j-1}})\nonumber\\
&-(n-j+2)\sum^{j-1}_{i=1}\frac{k_1}{k_4}|\xi_i|^4-\sum^{j-1}_{i=1}\frac{k_2}{k_4}|\xi_i|^{4+k_3}\nonumber\\
&+\lceil\xi_{j}\rfloor^{\frac{4 r-\kappa-r_{j}}{r}}h_j\lceil x_{j+1}\rfloor^{q_j}
+\lceil\xi_{j}\rfloor^{\frac{4 r-\kappa-r_{j}}{r}}f_j(\bar{x}_j)\nonumber\\
&+\sum^{j-1}_{i=1}\frac{\partial W_j}{\partial x_i}(h_i\lceil x_{i+1}\rfloor^{q_i}+f_i(\bar{x}_i))\nonumber\\
&+\frac{1}{2}\sum^j_{i,m=1}\left|\frac{\partial^2 W_j}{\partial x_i\partial x_m}\right|\left|g_i(\bar{x}_i)'g_m(\bar{x}_m)\right|.
\end{align}
Based on Propositions
\ref{predefined-time-2020-pro2}-\ref{predefined-time-2020-pro5}  in
Appendix, we estimate some terms in (\ref{predefined-time-2020-2})
and obtain that
{\small
\begin{align*}
&\mathcal{L}_jV_{j}(\bar{x}_{j})\\
\leq&(\alpha_{1,j}+\alpha_{2,j}+\alpha_{3,j}+\alpha_{4,j})|\xi_j|^4-(n-j+1)\sum^{j-1}_{i=1}\frac{k_1}{k_4}|\xi_i|^4\\
&-\sum^{j-1}_{i=1}\frac{k_2}{k_4}|\xi_i|^{4+k_3}+\lceil\xi_{j}\rfloor^{\frac{4 r-\kappa-r_{j}}{r}}h_j\lceil x_{j+1}^*\rfloor^{q_j}\\
&+\lceil\xi_{j}\rfloor^{\frac{4 r-\kappa-r_{j}}{r}}h_j( \lceil x_{j+1}\rfloor^{q_j}-\lceil x_{j+1}^*\rfloor^{q_j}).
\end{align*}
}
It is easy to see that the virtual controller
$$
x_{j+1}^*=-\beta_j\lceil\xi_{j}\rfloor^{\frac{r_{j+1}}{r}}
$$
and
\begin{align*}
\beta_j=&\left(\frac{(n-j+1)k_1}{\overline{h}_jk_4}+\frac{\alpha_{1,j}+\alpha_{2,j}+\alpha_{3,j}+\alpha_{4,j}}{\overline{h}_j}\right.\\
&\left.+\frac{k_2}{k_4\overline{h}_j}|\xi_j|^{k_3}\right)^{\frac{1}{q_j}}.
\end{align*}
Then
\begin{align}
&\mathcal{L}_jV_{j}(\bar{x}_{j})\nonumber\\
\leq&-(n-j+1)\sum^{j}_{i=1}\frac{k_1}{k_4}|\xi_i|^4-\sum^{j}_{i=1}\frac{k_2}{k_4}|\xi_i|^{4+k_3}\nonumber\\
&+\lceil\xi_{j}\rfloor^{\frac{4 r-\kappa-r_{j}}{r}}h_j(\lceil x_{j+1}\rfloor^{q_j}-\lceil x_{j+1}^*\rfloor^{q_j}).
\label{eqnjnj}
\end{align}
From Steps 1-3, we have proved that  (\ref{eqnjnj}) hold for any $j\in \{1,2,\cdots, n\}$.
Hence, at the $n$th step, one concludes that
\begin{align*}
&\mathcal{L}_nV_{n}(\bar{x}_{n})\\
\leq&-\sum^{n-1}_{i=1}\frac{k_1}{k_4}|\xi_i|^4-\sum^{n-1}_{i=1}\frac{k_2}{k_4}|\xi_i|^{4+k_3}+\lceil\xi_{n}\rfloor^{\frac{4 r-\kappa-r_n}{r}}h_n\lceil  u\rfloor^{q_n}  \nonumber\\
&+(\alpha_{1,n}+\alpha_{2,n}+\alpha_{3,n}+\alpha_{4,n})|\xi_n|^4 \nonumber\\
\leq&-\sum^{n}_{i=1}\frac{k_1}{k_4}|\xi_i|^4-\sum^{n}_{i=1}\frac{k_2}{k_4}|\xi_i|^{4+k_3}
\end{align*}
with the Lyapunov function as
$$
V_n(\bar{x}_n)=V_{n-1}(\bar{x}_{n-1})+W_n(\bar{x}_n)
$$
and
$$
W_n(\bar{x}_n)=\int^{x_n}_{x_n^*}\lceil\lceil s \rfloor^{\frac{r}{r_j}}-\lceil x_n^{*}\rfloor^{\frac{r}{r_n}}\rfloor^{\frac{4 r-\kappa-r_n}{r}}\,ds.
$$
Consequently, the control $u$ can be designed as
{\small
\begin{eqnarray*}
u&=&-\beta_n \lceil\xi_{n}\rfloor^{\frac{r_{n+1}}{r}}\\
&=&-\left(\frac{k_1}{\bar{h}_nk_4}+\frac{k_2}{k_4\bar{h}_n}|\xi_n|^{k_3}+\frac{\alpha_{1,n}+\alpha_{2,n}+\alpha_{3,n}+\alpha_{4,n}}
{\bar{h}_n}\right)^{\frac{1}{q_n}}\\
&&\cdot \lceil\xi_{n}\rfloor^{\frac{r_{n+1}}{r}}.
\end{eqnarray*}}
Based on the proof of Proposition \ref{predefined-time-2020-pro4}, we can get that
\begin{align*}
    \int^{x_j}_{x_j^*}\lceil\lceil s \rfloor^{\frac{r}{r_j}}-\lceil x_j^{*}\rfloor^{\frac{r}{r_j}}\rfloor^{\frac{4 r-\kappa-r_j}{r}}ds
    \leq2^{1-\frac{r_j}{r}}|\xi_j|^{\frac{4 r-\kappa}{r}}
\end{align*}
and
\begin{align*}
    V_n\leq2\sum^n_{i=1}|\xi_i|^{\frac{4 r-\kappa}{r}}.
\end{align*}
Setting $k_3\geq-\frac{\kappa}{r}$  and  using  Lemma \ref{predefined-time-2020-lem5} arrive at
\begin{align*}
&-\sum^{n}_{i=1}|\xi_i|^{4+k_3}\leq -n^{\frac{-rk_3-\kappa}{4r-\kappa}}(\sum^{n}_{i=1}|\xi_i|^{\frac{4r-\kappa}{r}})^{\frac{r(4+k_3)}{4r-\kappa}}\\
\leq&-n^{\frac{-rk_3-\kappa}{4r-\kappa}}\frac{1}{2}^{\frac{r(4+k_3)}{4r-\kappa}}V_n^{\frac{r(4+k_3)}{4r-\kappa}}
\end{align*}
and
$$
-\sum^{n}_{i=1}|\xi_i|^4\leq -\frac{1}{2}^{\frac{4r}{4r-\kappa}}V_n^{\frac{4r}{4r-\kappa}}.
$$
Accordingly, one has
\begin{align*}
    \mathcal{L}_nV_{n}(\bar{x}_{n})
    \leq&-\frac{k_1}{k_4}\left(\frac{V_n}{2}\right)^{\frac{4r}{4r-\kappa}}
    -\frac{k_2}{k_4}n^{\frac{-rk_3-\kappa}{4r-\kappa}}\left(\frac{V_n}{2}\right)^{\frac{4r+rk_3}{4r-\kappa}}.
\end{align*}
Set $k_3>-\frac{\kappa}{r}$, $b_1=\frac{4r}{4r-\kappa}$, $b_2=\frac{4r+rk_3}{4r-\kappa}$, $a=2^{-b1}k_1$, $k_1>\frac{2^{b1}}{1-b_1}$, $k_2=\frac{a-ab_1}{(a-ab_1-1)(b_2-1)}n^{\frac{rk_3+\kappa}{4r-\kappa}}2^{b_2}$ . The control parameter $k_4$ can be chosen as an arbitrary positive number. Then, considering Corollary \ref{predefined-time-cor3},
the origin of the closed-loop system is stochastically  predefined-time stabilizable and $\sup_{x_0\in\mathcal{R}^n}\mathcal{E}T(x_0)\leq k_4$.

\section{Simulation Examples}

In this section, we present two examples to illustrate the validity of our  main results.

\begin{example}
In this example, we consider the following one-dimensional system
\begin{align}\label{predefined-time-2020-ex1-system1}
dx=(x^{\frac{5}{3}}+u)dt+x^2dw.
\end{align}
For system (\ref{predefined-time-2020-ex1-system1}), the existing
design schemes \cite{2020TAC,Songzhibao} can not solve its
 stochastic predefined-time stabilization. Through the design method proposed in Corollary
\ref{predefined-time-cor3}, the controller can be constructed as
\begin{align*}
u=-x^{\frac{5}{3}}-\frac{1}{2}x^3-\frac{k_1}{2k_4}x^{\frac{1}{3}}-\frac{k_2}{2k_4}x^{k_3},
\end{align*}
where $k_1>4$, $k_2=\frac{k_1}{3(\frac{k_1}{4}-1)\frac{k_3-1}{2}}$, $k_3>1$, and $k_4$ is a positive control parameter that can be adjusted arbitrarily. Then, the mathematical expectation of the first time to reach the equilibrium point must be less than $k_4$. Here we choose $x(0)=1$, $k_1=4.1$, $k_3=3$. In order to compare the convergence speed, $k_4$ is selected as a variety of different values of $4$, $2$ and $0.5$. The corresponding trajectories of the states for  these groups of control experiments are all described in Figure \ref{example1}.
\begin{figure}[!htb]
  \centering
  \includegraphics[width=3in]{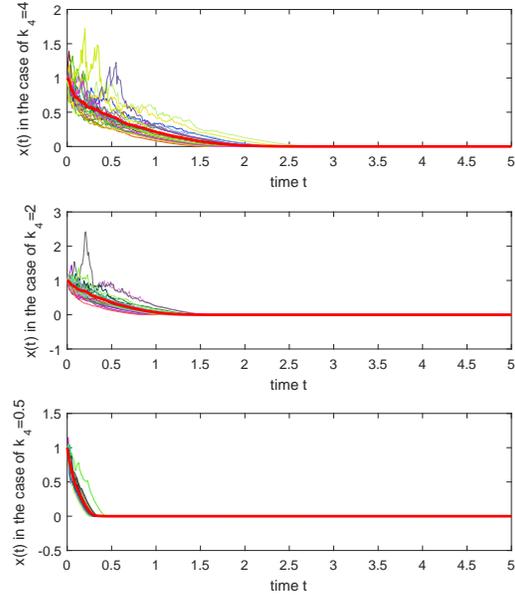}
  \caption{ The state $x(t)$ of the closed-loop system (\ref{predefined-time-2020-ex1-system1}) in the cases of $k_4=4,2,0.5$.}
  \label{example1}
\end{figure}
 We have completed $25$ random experiments for  $k_4=4,2,0.5$, respectively.
 From Figure \ref{example1}, we can find  that  smaller $k_4$   leads to  faster convergence speed. It is obvious that the rapidity of the convergence time meets the requirement of the stochastic fixed-time stabilization proposed in this paper.
\end{example}

\begin{example}
\begin{figure}[!htb]
  \centering
  \includegraphics[width=3in]{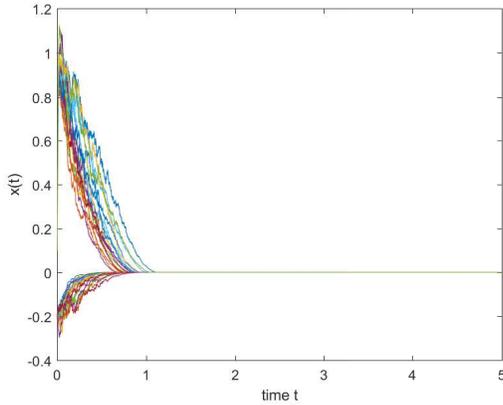}
  \caption{$x(t)$ of the closed-loop system (\ref{predefined-time-2020-ex1-system2}) in the cases of $k_4=1$.}
  \label{example2}
\end{figure}
A simulation example is given to illustrate how to apply Theorem \ref{predefined-time-th2}. Consider the following
nonlinear stochastic system:
\begin{align}\label{predefined-time-2020-ex1-system2}
\begin{cases}
dx_1=(h_1\lceil x_2 \rfloor^{\frac{5}{3}}-\lceil x_1 \rfloor^{\frac{3}{4}})dt+\sin(x_1)|x_1|dw,\\
dx_2=h_2\lceil u \rfloor^{\frac{4}{3}}\,dt.
\end{cases}
\end{align}
Clearly, Assumptions \ref{predefined-time-2020-ass1} and \ref{predefined-time-2020-ass2} are  satisfied with $\underline{h}_2=1$, $\overline{h}_2=2$, $\phi_1=\varphi_1=1$, $\kappa=-\frac{1}{4}$, $\varpi_{11}=0$, $\rho_{11}=\frac{1}{8}$. According to Theorem \ref{predefined-time-th2}, we can select the control parameters as $k_1=65.6$, $k_2=494.6$, $k_3=3.1$, $k_4=3$, and the control input can be designed as $u=-(10.9+3.5\beta_1(x_1)+82.4|\xi_2|^{0.75})\lfloor\xi_2\rceil^{\frac{1}{8}}$ with $\beta_1(x_1)=(47.1+165|\xi_1|^{3.1})^{0.6}$. Figure \ref{example2} shows  that the designed controller renders the origin of system (\ref{predefined-time-2020-ex1-system2}) stochastic predefined-time stabilization with $\mathcal{E}T_f\leq k_4$. This demonstrates the usability of the given design method.
\end{example}

\section{Conclusion}
In this paper, we have obtained  several Lyapunov-type results about
stochastic predefined-time stabilization of stochastic nonlinear
systems, which extend the results of finite-time and fixed-time
stabilization. Furthermore, we have presented   feasible conditions
and provided a constructive solution to stochastic predefined-time
stabilization of a class of stochastic nonlinear systems in
strict-feedback form. Two examples have been given to demonstrate
the validness of the obtained results.

\section{Appendix}

The following proposition can be obtained from \cite{Songzhibao}.
\begin{proposition}\label{predefined-time-2020-pro1}
The following statements accordingly hold when the positive definite and proper Lyapunov function $V_j$ and $W_j$ are defined by (\ref{predefined-time-2020-vv1}) and (\ref{predefined-time-2020-vv2}) respectively.
{\small\begin{align*}
&\frac{\partial W_j}{\partial x_j}=\lceil\xi_{j}\rfloor^{\frac{4 r-\kappa-r_{j}}{r}},\\
&\frac{\partial W_j}{\partial x_i}=-\frac{4r-\kappa-r_j}{r}\int^{x_j}_{x_j^*}|\lceil s \rfloor^{\frac{r}{r_j}}-\lceil x_j^{*}\rfloor^{\frac{r}{r_j}}|^{\frac{3 r-\kappa-r_j}{r}}ds\\
&\ \ \cdot \frac{\partial \lceil x_{j}^{*}\rfloor^{\frac{r}{r_j}}}{\partial x_i},\\
&\frac{\partial^2 W_j}{\partial^2 x_j}=\frac{4 r-\kappa-r_{j}}{r_j}|x_j|^{\frac{r-r_j}{r_j}}|\xi_{j}|^{\frac{3 r-\kappa-r_{j}}{r}},\\
&\frac{\partial^2 W_j}{\partial x_i\partial x_j}=\frac{\partial^2 W_j}{\partial x_j\partial x_i}=-\frac{4r-\kappa-r_j}{r}|\xi_j|^{\frac{3 r-\kappa-r_j}{r}}\frac{\partial \lceil x_{j}^{*}\rfloor^{\frac{r}{r_j}}}{\partial x_i},\\
&\frac{\partial^2 W_j}{\partial x_i \partial x_m}=\frac{(4r-\kappa-r_j)(3r-\kappa-r_j)}{r^2}\frac{\partial \lceil x_{j}^{*}\rfloor^{\frac{r}{r_j}}}{\partial x_i}\frac{\partial \lceil x_{j}^{*}\rfloor^{\frac{r}{r_j}}}{\partial x_m}\\
&\cdot \int^{x_j}_{x_j^*}\lceil\lceil s \rfloor^{\frac{r}{r_j}}-\lceil x_j^{*}\rfloor^{\frac{r}{r_j}}\rfloor^{\frac{2 r-\kappa-r_j}{r}}\,ds-\frac{4r-\kappa-r_j}{r}\frac{\partial^2 \lceil x_{j}^{*}\rfloor^{\frac{r}{r_j}}}{\partial x_i\partial x_m}\\
&\cdot \int^{x_j}_{x_j^*}|\lceil s \rfloor^{\frac{r}{r_j}}-\lceil x_j^{*}\rfloor^{\frac{r}{r_j}}|^{\frac{3 r-\kappa-r_j}{r}}\,ds,
\end{align*}}
where $i,m<j$.
\end{proposition}

Next, we need to prove some necessary propositions listed below   in order  to estimate some terms in (\ref{predefined-time-2020-2}).

\begin{proposition}\label{predefined-time-2020-pro2}
There exists a non-negative smooth function $\alpha_{1j}(\bar{x}_j)$ such that
\begin{align*}
    &\lceil\xi_{j-1}\rfloor^{\frac{4 r-\kappa-r_{j-1}}{r}}h_{j-1}(\lceil x_j\rfloor^{q_{j-1}}-\lceil x_j^*\rfloor^{q_{j-1}})\\
    \leq&
    \frac{k_1}{4k_4}|\xi_{j-1}|^4+\alpha_{1j}(\bar{x}_j)|\xi_j|^4.
\end{align*}
\end{proposition}

\begin{proposition}\label{predefined-time-2020-pro3}
There exists a non-negative smooth function $\alpha_{2j}(\bar{x}_j)$ such that
\begin{align*}
    \lceil\xi_{j}\rfloor^{\frac{4 r-\kappa-r_{j}}{r}}f_j(\bar{x}_j)
    \leq\frac{k_1}{4k_4}\sum^{j-1}_{i=1}|\xi_i|^{4}+\alpha_{2j}(\bar{x}_j)|\xi_j|^4.
\end{align*}
\end{proposition}

\begin{proposition}\label{predefined-time-2020-pro4}
There exists a non-negative smooth function $\alpha_{3j}(\bar{x}_j)$ such that
\small{\begin{equation*}
\sum^{j-1}_{i=1}\frac{\partial W_j}{\partial x_i}(h_i\lceil x_{i+1}\rfloor^{q_i}+f_i(\bar{x}_i))
\leq\frac{k_1}{4k_4}\sum^{j-1}_{i=1}|\xi_i|^4+\alpha_{3j}(\bar{x}_j)|\xi_j|^4.
\end{equation*}}
\end{proposition}

\begin{proposition}\label{predefined-time-2020-pro5}
There exists a non-negative smooth function $\alpha_{4j}(\bar{x}_j)$ such that
\small{\begin{equation*}
\frac{1}{2}\sum^j_{i,m=1}\left|\frac{\partial^2 W_j}{\partial x_i\partial x_m}\right|\left|g_i(\bar{x}_i)'g_m(\bar{x}_m)\right|\leq
\frac{k_1}{4k_4}\sum^{j-1}_{i=1}|\xi_i|^4+\alpha_{4,j}|\xi_j|^4.
\end{equation*}}
\end{proposition}

{\bf Proof of Proposition \ref{predefined-time-2020-pro2}}
\begin{small}\begin{align*}
    &\lceil\xi_{j-1}\rfloor^{\frac{4 r-\kappa-r_{j-1}}{r}}h_{j-1}(\lceil x_j\rfloor^{q_{j-1}}-\lceil x_j^*\rfloor^{q_{j-1}})\nonumber\\
    \leq&|\xi_{j-1}|^{\frac{4 r-\kappa-r_{j-1}}{r}}\bar{h}_{j-1}|(\lceil x_{j}\rfloor^{q_{j-1}\frac{r}{r_j}})^{\frac{r_j}{r}}-(\lceil x_{j}^*\rfloor^{q_{j-1}\frac{r}{r_j}})^{\frac{r_j}{r}}|.
\end{align*}
\end{small}
Based on Lemmas \ref{predefined-time-2020-lem2} and \ref{predefined-time-2020-lem3}, and $q_{j-1}r_j=r_{j-1}+\kappa$, the above inequality leads to
\begin{small}\begin{align*}
&\lceil\xi_{j-1}\rfloor^{\frac{4 r-\kappa-r_{j-1}}{r}}\bar{h}_{j-1}(\lceil x_j\rfloor^{q_{j-1}}-\lceil x_j^*\rfloor^{q_{j-1}})\nonumber\\
    \leq&  |\xi_{j-1}|^{\frac{4 r-\kappa-r_{j-1}}{r}}\bar{h}_{j-1}2^{1-\frac{r_{j-1}+\kappa}{r}}
    |\lceil x_{j}\rfloor^{\frac{r}{r_j}}-\lceil x_{j}^{*}\rfloor^{\frac{r}{r_j}}|^{\frac{r_{j-1}+\kappa}{r}}\\
    \leq&\frac{k_1}{4k_4}|\xi_{j-1}|^4+\bar{h}_{j-1}2^{1-\frac{r_{j-1}+\kappa}{r}}\frac{r_{j-1}+\kappa}{4r}\Phi_j^{-\frac{4 r-\kappa-r_{j-1}}{r_{j-1}+\kappa}}|\xi_j|^4,
\end{align*}\end{small}
where the real-valued function $\Phi_j$ can be selected as
$\Phi_j=\frac{rk_1}{k_4(4
r-\kappa-r_{j-1})\bar{h}_{j-1}2^{1-\frac{r_{j-1}+\kappa}{r}}}$. The
proof is completed.

{\bf Proof of Proposition \ref{predefined-time-2020-pro3}}
In view of Assumption \ref{predefined-time-2020-ass1}, we have
\begin{align*}
    &\lceil\xi_{j}\rfloor^{\frac{4 r-\kappa-r_{j}}{r}}f_j(\bar{x}_j)\\
    \leq&|\xi_{j}|^{\frac{4 r-\kappa-r_{j}}{r}}\bar{\phi}_j(\bar{x}_j)\sum^j_{i=1}||\xi_i|+|\beta_{i-1}|^{\frac{r}{r_i}}|\xi_{i-1}||^{\frac{r_j+\kappa}{r}}
\end{align*}
with $\xi_i=0$. Note that $r_j+\kappa<r$. It is obtained from Lemma \ref{predefined-time-2020-lem4} that
\begin{small}\begin{align*}
    &\lceil\xi_{j}\rfloor^{\frac{4 r-\kappa-r_{j}}{r}}f_j(\bar{x}_j)\\
    \leq&|\xi_{j}|^{\frac{4 r-\kappa-r_{j}}{r}}\bar{\phi}_j(\bar{x}_j)\left(\sum^{j-1}_{i=1}(1+|\beta_{i}|^{\frac{r_j+\kappa}{r_{i+1}}})|\xi_i|^{\frac{r_j+\kappa}{r}}+|\xi_j|^{\frac{r_j+\kappa}{r}}\right).
\end{align*}\end{small}
From Lemma \ref{predefined-time-2020-lem3},
\begin{align*}
    \lceil\xi_{j}\rfloor^{\frac{4 r-\kappa-r_{j}}{r}}f_j(\bar{x}_j)
    \leq\frac{k_1}{4k_4}\sum^{j-1}_{i=1}|\xi_i|^4
    +\alpha_{2j}(\bar{x}_j)|\xi_j|^{4}
   \end{align*}
where $\alpha_{2j}$ is a non-negative smooth function.
The proof is completed.

{\bf Proof of Proposition \ref{predefined-time-2020-pro4}}

According to Lemma \ref{predefined-time-2020-lem4} and the proof of Proposition \ref{predefined-time-2020-pro3}, we have
$f_i(\bar{x}_i)\leq\bar{\phi}_i(\bar{x}_i)\left(\sum^{i-1}_{m=1}(1+|\beta_{m}|^{\frac{r_i+\kappa}{r_{m+1}}})|\xi_m|^{\frac{r_i+\kappa}{r}}+|\xi_i|^{\frac{r_i+\kappa}{r}}\right)$,
and
$
\lceil x_{i+1}\rfloor^{q_i}
\leq|\xi_{i+1}|^{\frac{r_i+\kappa}{r}}+|\beta_i|^{\frac{r_i+\kappa}{r_{i+1}}}|\xi_i|^{\frac{r_i+\kappa}{r}}.
$
Considering Proposition \ref{predefined-time-2020-pro1}, we can get that
{\small\begin{align}\label{predefined-time-2020-nn3}
&\sum^{j-1}_{i=1}\frac{\partial W_j}{\partial x_i}(h_i\lceil x_{i+1}\rfloor^{q_i}+f_i(\bar{x}_i))\nonumber\\
\leq&\sum^{j-1}_{i=1}\left|\frac{4r-\kappa-r_j}{r}\frac{\partial \lceil x_{j}^{*}\rfloor^{\frac{r}{r_j}}}{\partial x_i}\int^{x_j}_{x_j^*}|\lceil s \rfloor^{\frac{r}{r_j}}-\lceil x_j^{*}\rfloor^{\frac{r}{r_j}}|^{\frac{3 r-\kappa-r_j}{r}}ds\right|\nonumber\\
&\Bigg[\bar{\phi}_i(\bar{x}_i)\left(\sum^{i-1}_{m=1}(1+|\beta_{m}|^{\frac{r_i+\kappa}{r_{m+1}}})|\xi_m|^{\frac{r_i+\kappa}{r}}+|\xi_i|^{\frac{r_i+\kappa}{r}}\right)\nonumber\\
&+|\xi_{i+1}|^{\frac{r_i+\kappa}{r}}+|\beta_i|^{\frac{r_i+\kappa}{r_{i+1}}}|\xi_i|^{\frac{r_i+\kappa}{r}}\Bigg].
\end{align}}
Using Lemma \ref{predefined-time-2020-lem2}, one has
\begin{align}\label{predefined-time-2020-nn4}
&\int^{x_j}_{x_j^*}\lceil\lceil s \rfloor^{\frac{r}{r_j}}-\lceil x_j^{*}\rfloor^{\frac{r}{r_j}}\rfloor^{\frac{3 r-\kappa-r_j}{r}}ds
\leq2^{1-\frac{r_j}{r}}|\xi_j|^{\frac{3 r-\kappa}{r}}.
\end{align}
Next, by  induction, we can estimate
\begin{align}\label{predefined-time-2020-nn5}
\left|\frac{\partial \lceil x_{j}^{*}\rfloor^{\frac{r}{r_j}}}{\partial x_i}\right|\leq \Theta_{1ji}\sum^{j-1}_{l=1}|\xi_l|^{\frac{r-r_i}{r}},\ j>i\geq1,
\end{align}
where $\Theta_{1ji}$ is a non-negative smooth function.
According to (\ref{predefined-time-2020-nn3}), (\ref{predefined-time-2020-nn4}) and (\ref{predefined-time-2020-nn5}), we have
\begin{align*}
&\sum^{j-1}_{i=1}\frac{\partial W_j}{\partial x_i}(h_i\lceil x_{i+1}\rfloor^{q_i}+f_i(\bar{x}_i))\nonumber\\
\leq&\sum^{j-1}_{i=1}\frac{4r-\kappa-r_j}{r}2^{1-\frac{r_j}{r}}|\xi_j|^{\frac{3r-\kappa}{r}}\Theta_{1ji}\sum^{j-1}_{l=1}|\xi_l|^{\frac{r-r_i}{r}}
\Theta_{2i}\nonumber\\
&\sum^{j}_{m=1}|\xi_m|^{\frac{r_i+\kappa}{r}}\nonumber\\
\leq&\frac{k_1}{4k_4}\sum^{j-1}_{i=1}|\xi_i|^4+\alpha_{3j}(\bar{x}_j)|\xi_j|^4.
\end{align*}
where $\alpha_{3j}$ is a non-negative smooth function.
Proposition \ref{predefined-time-2020-pro4} is proved.

{\bf Proof of Proposition \ref{predefined-time-2020-pro5}}
{\small
\begin{align*}
&\frac{1}{2}\sum^j_{i,m=1}\left|\frac{\partial^2 W_j}{\partial x_i\partial x_m}\right|\left|g_i(\bar{x}_i)'g_m(\bar{x}_m)\right|\\
=&\frac{1}{2}\left|\frac{\partial^2 W_j}{\partial^2 x_j}\right|\left|g_j(\bar{x}_j)'g_j(\bar{x}_j)\right|
+\sum^{j-1}_{i=1}\left|\frac{\partial^2 W_j}{\partial x_j\partial x_i}\right|\left|g_j(\bar{x}_j)'g_i(\bar{x}_i)\right|\\
&+\frac{1}{2}\sum^{j-1}_{i,m=1}\left|\frac{\partial^2 W_j}{\partial x_i\partial x_m}\right|\left|g_i(\bar{x}_i)'g_m(\bar{x}_m)\right|.
\end{align*}}
Firstly, based on Propositions \ref{predefined-time-2020-pro1}, \ref{predefined-time-2020-pro3} and  Lemma \ref{predefined-time-2020-lem5}, we have
\small{\begin{align}\label{predefined-time-2020-nn15}
&\frac{1}{2}\left|\frac{\partial^2 W_j}{\partial^2 x_j}\right|\left|g_j(\bar{x}_j)'g_j(\bar{x}_j)\right|\nonumber\\
\leq&\frac{(4r-\kappa-r_j)j}{2r_j}\Bigg(|\xi_j|^{\frac{r-r_j}{r}}+|\beta_{j-1}|^{\frac{r-r_j}{r_j}}|\xi_{j-1}|^{\frac{r-r_j}{r}}\Bigg)\bar{\varphi}_j(\bar{x}_j)^2
\nonumber\\&|\xi_j|^{\frac{3r-\kappa-r_j}{r}}
\Bigg(\sum^{j-1}_{s=1}(1+|\beta_s|^{\frac{2r_j+\kappa}{r_{s+1}}})|\xi_s|^{\frac{2r_j+\kappa}{r}}+|\xi_j|^{\frac{2r_j+\kappa}{r}}\Bigg)\nonumber\\
\leq&\Xi_{1j}(\zeta_{1j},\bar{x}_j)|\xi_j|^4+\zeta_{1j}\sum^{j-1}_{s=1}|\xi_s|^4,
\end{align}}
where $\Xi_{1j}(\zeta_{1j},\bar{x}_j)$ is a positive real function and $\zeta_{1j}$ can be an any positive real number.
Next, we will consider
\begin{align*}
&\sum^{j-1}_{i=1}\left|\frac{\partial^2 W_j}{\partial x_j\partial x_i}\right|\left|g_j(\bar{x}_j)'g_i(\bar{x}_i)\right|\\
\leq&\sum^j_{i=1}\frac{4r-\kappa-r_j}{r}|\xi_j|^{\frac{3 r-\kappa-r_j}{r}}\Theta_{1ji}\sum^{j-1}_{l=1}|\xi_l|^{\frac{r-r_i}{r}}\left|g_j(\bar{x}_j)'g_i(\bar{x}_i)\right|.
\end{align*}
Note that
$$\left|g_j(\bar{x}_j)'g_i(\bar{x}_i)\right|\leq
\bar{\varphi}_i(\bar{x}_i)\bar{\varphi}_j(\bar{x}_j)\Theta_{2ji}\sum^{j}_{s=1}|\xi_s|^{\frac{r_j+r_i+\kappa}{r}}$$
with $\Theta_{2ji}$ being a positive  continuous function.
Summarizing the above analysis leads to that
\begin{align}\label{predefined-time-2020-nn18}
&\sum^{j-1}_{i=1}\left|\frac{\partial^2 W_j}{\partial x_j\partial x_i}\right|\left|g_j(\bar{x}_j)'g_i(\bar{x}_i)\right|\nonumber\\
\leq&\sum^j_{i=1}\frac{(4r-\kappa-r_j)j}{r}|\xi_j|^{\frac{3r-\kappa-r_j}{r}}\Theta_{1ji}\Theta_{2ji}\bar{\varphi}_i(\bar{x}_i)\bar{\varphi}_j(\bar{x}_j)\nonumber\\
&\sum^{j-1}_{l=1}|\xi_l|^{\frac{r-r_i}{r}}\sum^{j}_{s=1}|\xi_s|^{\frac{r_j+r_i+\kappa}{r}}\nonumber\\
\leq&\Xi_{2j}(\zeta_{2j},\bar{x}_j)|\xi_j|^4+\zeta_{2j}\sum^{j-1}_{l=1}|\xi_l|^4,
\end{align}
where $\Xi_{2j}(\zeta_{2j},\bar{x}_j)$ is a positive real function and $\zeta_{2j}$ can be an any positive real number.

Finally, we have the following estimation:
\begin{align*}
&\frac{1}{2}\sum^{j-1}_{i,m=1}\left|\frac{\partial^2 W_j}{\partial x_i\partial x_m}\right|\left|g_i(\bar{x}_i)'g_m(\bar{x}_m)\right|\\
\leq&\sum^{j-1}_{i,m=1}\Bigg|\frac{(4r-\kappa-r_j)(3r-\kappa-r_j)}{r^2}\frac{\partial \lceil x_{j}^{*}\rfloor^{\frac{r}{r_j}}}{\partial x_i}\frac{\partial \lceil x_{j}^{*}\rfloor^{\frac{r}{r_j}}}{\partial x_m}\\
&\int^{x_j}_{x_j^*}\lceil\lceil s \rfloor^{\frac{r}{r_j}}-\lceil x_j^{*}\rfloor^{\frac{r}{r_j}}\rfloor^{\frac{2 r-\kappa-r_j}{r}}ds-\frac{4r-\kappa-r_j}{r}\frac{\partial^2 \lceil x_{j}^{*}\rfloor^{\frac{r}{r_j}}}{\partial x_i\partial x_m}\\
&\int^{x_j}_{x_j^*}|\lceil s \rfloor^{\frac{r}{r_j}}-\lceil x_j^{*}\rfloor^{\frac{r}{r_j}}|^{\frac{3 r-\kappa-r_j}{r}}ds\Bigg|\bar{\varphi}_i(\bar{x}_i)\bar{\varphi}_m(\bar{x}_m)\\
&\Theta_{2im}\sum^{j-1}_{s=1}|\xi_s|^{\frac{r_m+r_i+\kappa}{r}}.
\end{align*}
According to (\ref{predefined-time-2020-nn4}) and (\ref{predefined-time-2020-nn5}), we have
$\int^{x_j}_{x_j^*}\lceil\lceil s \rfloor^{\frac{r}{r_j}}-\lceil x_j^{*}\rfloor^{\frac{r}{r_j}}\rfloor^{\frac{2 r-\kappa-r_j}{r}}ds\leq2^{1-\frac{r_j}{r}}|\xi_j|^{\frac{2 r-\kappa}{r}}
$
and
$\left|\frac{\partial \lceil x_{j}^{*}\rfloor^{\frac{r}{r_j}}}{\partial x_i}\right|\left|\frac{\partial \lceil x_{j}^{*}\rfloor^{\frac{r}{r_j}}}{\partial x_m}\right|
\leq\Theta_{3jim}\sum^{j-1}_{l=1}|\xi_l|^{\frac{2r-r_i-r_m}{r}}$
with $\Theta_{3jim}$ being a positive real function.
Therefore, we can get that
\begin{align}\label{predefined-time-2020-nn19}
&\sum^{j-1}_{i,m=1}\frac{(4r-\kappa-r_j)(3r-\kappa-r_j)}{r^2}\frac{\partial \lceil x_{j}^{*}\rfloor^{\frac{r}{r_j}}}{\partial x_i}\frac{\partial \lceil x_{j}^{*}\rfloor^{\frac{r}{r_j}}}{\partial x_m}\nonumber\\
&\int^{x_j}_{x_j^*}\lceil\lceil s \rfloor^{\frac{r}{r_j}}-\lceil x_j^{*}\rfloor^{\frac{r}{r_j}}\rfloor^{\frac{2 r-\kappa-r_j}{r}}ds\bar{\varphi}_i(\bar{x}_i)\bar{\varphi}_m(\bar{x}_m)\nonumber\\
&\Theta_{2im}\sum^{j-1}_{s=1}|\xi_s|^{\frac{r_m+r_i+\kappa}{r}}\nonumber\\
\leq&\Xi_{3j}(\zeta_{3j},\bar{x}_j)|\xi_j|^{4}+\zeta_{3j}\sum^{j-1}_{l=1}|\xi_l|^{4},
\end{align}
where  $\Xi_{3j}(\zeta_{3j},\bar{x}_j)$ is a positive real function and $\zeta_{3j}$ can be an any positive real number.
Now, we are in a position to consider $\sum^{j-1}_{i,m=1}\Bigg|\frac{4r-\kappa-r_j}{r}\frac{\partial^2 \lceil x_{j}^{*}\rfloor^{\frac{r}{r_j}}}{\partial x_i\partial x_m}\int^{x_j}_{x_j^*}|\lceil s \rfloor^{\frac{r}{r_j}}-\lceil x_j^{*}\rfloor^{\frac{r}{r_j}}|^{\frac{3 r-\kappa-r_j}{r}}ds\Bigg|\bar{\varphi}_i(\bar{x}_i)\bar{\varphi}_m(\bar{x}_m)
\Theta_{3im}\sum^{j-1}_{s=1}|\xi_s|^{\frac{r_m+r_i+\kappa}{r}}$.
Through calculation, we can obtain $\frac{\partial^2 \lceil x_{j}^{*}\rfloor^{\frac{r}{r_j}}}{\partial x_i\partial x_m}
\leq\sum^{j-1}_{s=1}\Theta_{4jims}|\xi_s|^{\frac{r-r_i-r_m}{r}}$
with
$\Theta_{4jims}=\prod^{j-1}_{l=s+1}|\beta_l^{\frac{r}{r_{l+1}}}|\Bigg|\frac{\partial^2 \beta_{s}^{\frac{r}{r_{s+1}}}}{\partial x_i\partial x_m}\Bigg||\xi_s|^{\frac{r_i+r_m}{r}}+I_{i=m}\prod^{j-1}_{l=s}|\beta_l^{\frac{r}{r_{l+1}}}|.$
As a consequence,
{\small\begin{align}\label{predefined-time-2020-nn20}
&\sum^{j-1}_{i,m=1}\Bigg|\frac{4r-\kappa-r_j}{r}\frac{\partial^2 \lceil x_{j}^{*}\rfloor^{\frac{r}{r_j}}}{\partial x_i\partial x_m}
\int^{x_j}_{x_j^*}\lceil\lceil s \rfloor^{\frac{r}{r_j}}-\lceil x_j^{*}\rfloor^{\frac{r}{r_j}}\rfloor^{\frac{3 r-\kappa-r_j}{r}}ds\Bigg|\nonumber\\
&\bar{\varphi}_i(\bar{x}_i)\bar{\varphi}_m(\bar{x}_m)\sum^m_{s=1}|x_s|^{\frac{2r_m+\kappa}{2r_s}}\sum^i_{s=1}|x_s|^{\frac{2r_i+\kappa}{2r_s}}\nonumber\\
&\leq\sum^{j-1}_{i,m=1}\frac{4r-\kappa-r_j}{r}2^{1-\frac{r_j}{r}}|\xi_j|^{\frac{3 r-\kappa}{r}}\bar{\varphi}_i(\bar{x}_i)\bar{\varphi}_m(\bar{x}_m)\sum^{j-1}_{s=1}\Theta_{4jims}\nonumber\\
&|\xi_s|^{\frac{r-r_i-r_m}{r}}\Theta_{3im}\sum^{j-1}_{s=1}|\xi_s|^{\frac{r_m+r_i+\kappa}{r}}\nonumber\\
&\leq\Xi_{4j}(\zeta_{4j},\bar{x}_j)|\xi_j|^{4}+\zeta_{4j}\sum^{j-1}_{s=1}|\xi_s|^4.
\end{align}}
where $\Xi_{4j}(\zeta_{4j},\bar{x}_j)$ is a positive real function and $\zeta_{4j}$ can be an   any   positive real number.
Based on (\ref{predefined-time-2020-nn15}), (\ref{predefined-time-2020-nn18}), (\ref{predefined-time-2020-nn19}) and (\ref{predefined-time-2020-nn20}), we can find  a  suitable    $\zeta_{ij}$ such that  $\sum^{4}_{i=1}\zeta_{ij}=\frac{k_1}{4k_4}$, which
completes the proof of Proposition \ref{predefined-time-2020-pro5}.

\end{document}